# Distributed Finite Time Termination of Ratio Consensus for Averaging in the presence of Delays

Mangal Prakash*, Saurav Talukdar*, Sandeep Attree, Sourav Patel and Murti V. Salapaka

*Both authors have contributed equally

## Abstract

Distributed averaging of agent initial conditions is a well-studied problem in context of networked systems where coordination amongst the agents is of paramount importance. The asymptotic nature of convergence of distributed averaging protocols and presence of communication delays, however, makes it challenging to implement in practical settings. It is important that agents develop the ability to detect on their own when average of the initial conditions of the agents is achieved within some pre-specified tolerance and stop further computations to avoid overhead expenses. This article presents a distributed finite time stopping criterion for distributed averaging using ratio consensus on a fixed interconnection topology (captured by a directed or undirected graph). The practical utility of the algorithm has been illustrated through simulations as well as experiments on a communication network realized with Raspberry Pis.

## I. INTRODUCTION

The problem of consensus in networks of autonomous agents has received significant attention from scientific community. The act of agreement, also often called consensus, is a prerequisite for many distributed control problems such as decentralized control of micro grids, localization of sensor networks, control of vehicular formations and many similar problems in robotics [1]. Typically, each agent participating in a consensus protocol holds parameter(s) of interest as it's state and keeps updating it based on the state of its neighbors. In this article, we are interested in a special case of consensus referred to as average consensus, where the agents need to synchronize their states to the average of the initial conditions. However, achieving average consensus in real time applications is challenging as it is typically achieved asymptotically.

Specifically, this article deals with computing average consensus in a distributed manner in finite time which can be applicable for practical engineered problems such as power allocation/frequency regulation in micro grids [2]. Naturally such a requirement translates into terminating the average consensus computations at the individual agent level when the values held by the agents are within a pre-specified bound from the desired average value. Our discussion pertains to the discrete time setting as most of the hardware implementations of such algorithms invariably assume a discrete time model. One of the first few algorithms to address the issue of finite time convergence was reported in [3] where max-min consensus was used to determine the proximity of the agent value to the desired average value in finite time. Sundaram and Hadjicostis proposed a finite time algorithm which relies on computing the minimal polynomial of the weight matrix [4]. This approach becomes cumbersome when a decentralized implementation of the algorithm is sought. It turns out that decentralized computation of

S. Talukdar is with the Department of Mechanical Engineering, University of Minnesota, Minneapolis, USA sauravtalukdar@umn.edu

M. Prakash, S. Attree, S. Patel and M. V. Salapaka are with the Department of Electrical and Computer Engineering, University of Minnesota, Minneapolis, USA praka027@umn.edu, attree002@umn.edu, murtis@umn.edu

minimal polynomial by each agent is much more demanding in terms of computational complexity and memory requirements than when done in a centralized fashion. [2] presented an algorithm for distributed finite time termination of ratio consensus [5], [6], which leveraged the simplicity of max-min consensus approach of [3] to distributively apportion resources in a AC micro grid setting. Recently, [7] analyzed the algorithm presented in [2] for distributed finite time termination of ratio consensus. However, the distributed finite time termination algorithms discussed above, neglected the presence of inherent latencies on communication channels, which happens to be a fundamental characteristic of any realistic large scale network.

Consensus in the presence of delays is discussed in [8], [9], [10], where a graphical framework is used to tackle delays/latencies by modeling them as virtual agents in the network. [11] uses a similar approach for handling the latencies but proposes a novel ratio consensus based algorithm, where the average of the initial conditions is achieved in a distributed and asymptotic manner by taking the ratio of two updates running at each node simultaneously. The highlighting feature of [11] is that the weight matrix is column stochastic, enabling the weights to be chosen in a distributed manner. Finite time termination of the ratio based average consensus in presence of delays presented in [11] is described in [12]. It relies on computing the rank of a Hankel matrix at each iteration by each node which can be challenging for large scale networks, thereby imposing serious constraints on computational complexity and memory requirements. The advantages and simplicity of Max-Min consensus was recently utilized in [13] for distributed finite-time termination of consensus in presence of delays as presented in [9]. This article extends max-min consensus for distributed finite time termination of ratio consensus in the presence of delays.

The rest of the article is organized as follows. In Section II, the basic definitions needed for subsequent development are presented. The setup for distributed averaging using ratio consensus is presented in Section III. In Section IV, the analytical results for distributed finite time termination of ratio consensus using max-min consensus is discussed followed by simulation results in Section V. Finally, concluding remarks are presented in Section VI.

## II. DEFINITIONS

In this section we present basic notions of graph theory and linear algebra which are essential for the subsequent developments. Detailed description of graph theory and linear algebra notions are available in [14] and [15] respectively.

*Definition 1:* (Directed and Undirected Graph) A directed graph $G$ is a pair $(V, E)$ where $V$ is a set of vertices or nodes and $E$ is a set of edges, which are ordered subsets of two distinct elements of $V$. If an edge from $j \in V$ to $i \in V$ exists then it is denoted as $(i, j) \in E$. An undirected graph $G$ is a pair $(V, E)$ where $V$ is a set of vertices or nodes and $E$ is a set of edges such that for every pair of distinct nodes $i \in V$ and $j \in V$, if $(i, j) \in E$ then $(j, i) \in E$.

*Definition 2:* (Path) In an undirected graph, a path between nodes $i \in V$ and $j \in V$ is a sequence of distinct edges of $G$ of the form $(i, k_1), (k_1, k_2), ..., (k_{m-1}, k_m), (k_m, j)$ and $(j, k_m), (k_m, k_{m-1}), ..., (k_2, k_1), (k_1, i)$.

In a directed graph, a directed path from node $i$ to $j$ exists if there is a sequence of distinct directed edges of $G$ of the form $(k_1, i), (k_2, k_1), ..., (j, k_m)$.

*Definition 3:* (Connected and Strongly Connected Graph) An undirected graph is connected if it has a path between each pair of distinct nodes $i$ and $j$. A directed graph is strongly connected if it has a directed path between each pair of distinct nodes $i$ and $j$.

*Definition 4:* (In-neighbor of node) Given a directed graph $G = (V, E)$, a node $j \in V$ is said to be an in-neighbor of node $i \in V$ if $(i, j) \in E$. The set of in-neighbors of node $i \in V$ is denoted by $N_i := \{j : (i, j) \in E\}$.

*Definition 5:* (Column Stochastic Matrix) A real $n \times n$ matrix $A = [a_{ij}]$ is called a column stochastic matrix if $1 \geq a_{ij} \geq 0$ for $1 \leq i, j \leq n$ and $\sum_{i=1}^{n} a_{ij} = 1$ for $1 \leq j \leq n$.

*Definition 6 (Irreducible Matrix):* A $N \times N$ matrix $A$ is said to be irreducible if for any $i, j \in \{1, ..., N\}$, there exist $m \in \mathbb{N}$ such that $(\mathbf{A}^m)(i, j) > 0$, that is, it is possible to reach any state from any other state in finite hops.

*Definition 7:* (Primitive Matrix) A non negative matrix $A$ is primitive if it is irreducible and has only one eigenvalue of maximum modulus.

*Definition 8:* (Diameter of a Graph) The diameter of a graph is defined as the longest shortest path between any two nodes in the network and is denoted by $D$.

### III. AVERAGING USING RATIO CONSENSUS

In this section, the key result from [5], which enables reaching average consensus using a ratio of two states maintained by each node is summarized. Consider a directed graph $G = \{V, E\}$ of $N$ nodes. Each node $i \in V$ maintains two states at time $k$, denoted by $x_i(k)$ (referred as numerator state of node $i$) and $y_i(k)$ (referred as denominator state of node $i$). Node $i$ updates its state at the $(k+1)^{th}$ iteration according to the following policy:

$$x_i(k+1) = p_{ii} x_i(k) + \sum_{j \in N_i} p_{ij} \sum_{r=0}^{\bar{\tau}} x_j(k-r) I_{k-r,ij}(r), \quad (1)$$

$$I_{k,ij}(\tau) = \begin{cases} 1, & \text{if } \tau_{ij}(k) = \tau \\ 0, & \text{if } \tau_{ij}(k) \neq \tau, \end{cases}$$

and

$$y_i(k+1) = p_{ii} y_i(k) + \sum_{j \in N_i} p_{ij} \sum_{r=0}^{\bar{\tau}} y_j(k-r) I_{k-r,ij}(r), \quad (2)$$

$$I_{k,ij}(\tau) = \begin{cases} 1, & \text{if } \tau_{ij}(k) = \tau \\ 0, & \text{if } \tau_{ij}(k) \neq \tau, \end{cases}$$

Thus, the delay between the neighbors is time varying but uniformly bounded by $\bar{\tau}$.

*Theorem 3.1:* [11] Suppose the weight matrix $P$ with $P(i,j) = p_{ij}$ associated with the directed graph $G$ is primitive and column stochastic with $P(i,i) > 0$ for all $i \in V$. Let the initial conditions for the numerator states be $\{x_1(0), x_2(0), ..., x_N(0)\}$ and for the denominator states be $\{y_i(0) = 1\}_{i \in V}$. Then the ratio $\frac{x_i(k)}{y_i(k)}$ asymptotically converges to $\frac{\sum_{i=1}^{N} x_i(0)}{N}$ for all $i = 1, ..., N$, that is, the ratio of the numerator and denominator states converge to the average of the initial conditions (referred as ratio consensus).

*Proof:* See [11] for proof. ∎

## IV. DISTRIBUTED FINITE TIME TERMINATION

In this section, first results based on the update rules (1) and (2) are established followed by the definitions and convergence of Max-Min consensus algorithms. Subsequently, a finite-time termination criterion for average consensus is developed based on these results. Let us consider the maximum and minimum value of the ratio of consensus protocols given by (1) and (2) over all nodes within a time horizon $\bar{\tau}$ from any time instant $k$ be given by,

$$M(k) := \max_{\substack{j \in V \\ r = \{0,1,2,\ldots,\bar{\tau}\}}} \frac{x_j(k-r)}{y_j(k-r)}, y_j(k-r) \neq 0, j \in V \quad (3)$$

$$m(k) := \min_{\substack{j \in V \\ r = \{0,1,2,\ldots,\bar{\tau}\}}} \frac{x_j(k-r)}{y_j(k-r)}, y_j(k-r) \neq 0, j \in V \quad (4)$$

*Remark 1:* Given that $y_j(0) = 1$ for all $j \in V$ and $P$ is a non-negative matrix, $y_j(k-r) \neq 0$ for all $k \in \mathbb{N}$ and $r \in \{0, 1, 2, \ldots, \bar{\tau}\}$.

*Lemma 4.1:* Consider the update rule (1) and (2) with the assumptions of Theorem 3.1. Then for all time instants $k' \geq k$ and for all $i \in V$,

$$\frac{x_i(k')}{y_i(k')} \leq M(k) \quad (5)$$

and

$$\frac{x_i(k')}{y_i(k')} \geq m(k). \quad (6)$$

*Proof:* We prove the theorem using strong induction. By definition of $M(k)$, for $k' = k$, the proof is trivial. Suppose it is asserted that for $k' = k + l', l' = \{1, 2, \ldots, l\}$, $\frac{x_i(k+l')}{y_i(k+l')} \leq M(k)$ for all $i \in V$. Then for $k' = k + l + 1$,

$$\frac{x_i(k+l+1)}{y_i(k+l+1)} = \frac{p_{ii}x_i(k+l) + \sum_{j \in N_i} p_{ij} \sum_{r=0}^{\bar{\tau}} x_j(k+l-r) I_{k+l-r,ij}(r)}{p_{ii}y_i(k+l) + \sum_{j \in N_i} p_{ij} \sum_{r=0}^{\bar{\tau}} y_j(k+l-r) I_{k+l-r,ij}(r)}$$

$$= \frac{p_{ii}\frac{x_i(k+l)}{y_i(k+l)} + \frac{1}{y_i(k+l)} \sum_{j \in N_i} p_{ij} \sum_{r=0}^{\bar{\tau}} x_j(k+l-r) I_{k+l-r,ij}(r)}{p_{ii} + \frac{1}{y_i(k+l)} \sum_{j \in N_i} p_{ij} \sum_{r=0}^{\bar{\tau}} y_j(k+l-r) I_{k+l-r,ij}(r)}.$$

Consider the case when $l > \bar{\tau}$. Then $k + l \geq k + l - r > k + \bar{\tau} - \bar{\tau} = k$. Using the inductive assertion, it follows that $\frac{x_j(k+l-r)}{y_j(k+l-r)} \leq M(k)$. This implies that $x_j(k+l-r) \leq M(k)y_j(k+l-r)$. Thus, it follows that

$$\frac{x_i(k+l+1)}{y_i(k+l+1)} \leq \frac{p_{ii}M(k) + \frac{1}{y_i(k+l)} \sum_{j \in N_i} p_{ij} \sum_{r=0}^{\bar{\tau}} y_j(k+l-r) I_{k+l-r,ij}(r) M(k)}{p_{ii} + \frac{1}{y_i(k+l)} \sum_{j \in N_i} p_{ij} \sum_{r=0}^{\bar{\tau}} y_j(k+l-r) I_{k+l-r,ij}(r)}.$$

. It follows that $\frac{x_i(k+l+1)}{y_i(k+l+1)} \leq M(k)$ for any $i \in V$.

Consider the case when $1 \leq l \leq \bar{\tau}$. Then $k + l - r \geq k + 1 - \bar{\tau} > k - \bar{\tau}$. Furthermore, suppose that $k \geq k + l - r \geq k + 1 - \bar{\tau} > k - \bar{\tau}$. Then from the definition of $M(k)$, it follows that $\frac{x_j(k+l-r)}{y_j(k+l-r)} \leq M(k)$. It follows that $x_j(k+l-r) \leq M(k)$. Now suppose $k + l - r > k$. Thus, it follows that $k + l \geq k + l - r > k$. Using the inductive assertion, it follows that

$\frac{x_j(k+l-r)}{y_j(k+l-r)} \leq M(k)$. This implies that $x_j(k+l-r) \leq M(k)y_j(k+l-r)$. It then follows that $\frac{x_i(k+l+1)}{y_i(k+l+1)} \leq M(k)$ for any $i \in V$. This proves (5). Similarly, (6) can be proven. This completes the proof. ∎

*Lemma 4.2:* Consider the update rule (1) and (2) with the assumptions of Theorem 3.1 with the initial time instant being $k$. Let $M(k) := \max_{j \in V} \frac{x_j(k)}{y_j(k)}$ and $m(k) := \min_{j \in V} \frac{x_j(k)}{y_j(k)}$. Let $i$ be a node such that $\frac{x_i(k')}{y_i(k')} < M(k)$ and let $j$ be a node such that $\frac{x_j(k')}{y_j(k')} > m(k)$ for some time instant $k' \geq k$. Then for all time instants $k'' \geq k'$, $\frac{x_i(k'')}{y_i(k'')} < M(k)$ and $\frac{x_j(k'')}{y_j(k'')} > m(k)$.

*Proof:* The proof is based on induction. It follows from the hypothesis that for time instant $k'' = k'$, $\frac{x_i(k')}{y_i(k')} < M(k)$. Suppose the assertion is that for some integer $q \geq 1$ and for $k'' = k' + q$, $\frac{x_i(k'+q)}{y_i(k'+q)} < M(k)$. This implies that for $k'' = k' + q$, $x_i(k'+q) < M(k)y_i(k'+q)$. Then for $k'' = k' + q + 1$,

$$\frac{x_i(k'+q+1)}{y_i(k'+q+1)} = \frac{p_{ii}x_i(k'+q) + \sum_{j \in N_i} p_{ij} \sum_{r=0}^{\bar{\tau}} x_j(k'+q-r) I_{k'+q-r,ij}(r)}{p_{ii}y_i(k'+q) + \sum_{j \in N_i} p_{ij} \sum_{r=0}^{\bar{\tau}} y_j(k'+q-r) I_{k'+q-r,ij}(r)} \quad (7)$$

Clearly, $k' - \bar{\tau} < k' + 1 - \bar{\tau} \leq k' + q - r \leq k' + q$. Since $k' \geq k$, it follows that $k' - \bar{\tau} \geq k - \bar{\tau}$. Thus, $k - \bar{\tau} \leq k' - \bar{\tau} < k' + 1 - \bar{\tau} \leq k' + q - r \leq k' + q$. Since $k' \geq k$. Now consider the case when $k - \bar{\tau} < k' + q - r \leq k$. Then using (3), it follows that $\frac{x_j(k'+q-r)}{x_j(k'+q-r)} \leq M(k)$ and so, $x_j(k'+q-r) \leq M(k)y_j(k'+q-r)$ for all $j \in V$. Using this fact together with the inductive assertion, it further follows from (7) that $\frac{x_i(k'+q+1)}{y_i(k'+q+1)} < \frac{p_{ii}M(k)y_i(k'+q) + \sum_{j \in N_i} p_{ij} \sum_{r=0}^{\bar{\tau}} M(k)y_j(k'+q-r) I_{k'+q-r,ij}(r)}{p_{ii}y_i(k'+q) + \sum_{j \in N_i} p_{ij} \sum_{r=0}^{\bar{\tau}} y_j(k'+q-r) I_{k'+q-r,ij}(r)}$. Thus, it follows that $\frac{x_i(k'+q+1)}{y_i(k'+q+1)} < M(k)$. Now consider the case when $k \leq k' + q - r \leq k' + q$. Then from *Lemma 4.1*, it follows that $\frac{x_j(k'+q-r)}{x_j(k'+q-r)} \leq M(k)$ and so, $x_j(k'+q-r) \leq M(k)y_j(k'+q-r)$ for all $j \in V$. Using this fact together with the inductive assertion, it further follows from (7) that $\frac{x_i(k'+q+1)}{y_i(k'+q+1)} < \frac{p_{ii}M(k)y_i(k'+q) + \sum_{j \in N_i} p_{ij} \sum_{r=0}^{\bar{\tau}} M(k)y_j(k'+q-r) I_{k'+q-r,ij}(r)}{p_{ii}y_i(k'+q) + \sum_{j \in N_i} p_{ij} \sum_{r=0}^{\bar{\tau}} y_j(k'+q-r) I_{k'+q-r,ij}(r)}$. Thus, it follows that $\frac{x_i(k'+q+1)}{y_i(k'+q+1)} < M(k)$.

The proof of the other inequality is similar to the proof above and is left to the reader. ∎

*Remark 2:* Note that $\frac{x_i(k')}{y_i(k')} \leq M(k)$, for $k' \geq k$ follows from Lemma 4.1. Lemma 4.2 emphasizes the strict inequality $\frac{x_i(k'')}{y_i(k'')} < M(k)$ for $k'' > k'$ if at any time instant $k'$, the inequality $\frac{x_i(k')}{y_i(k')} < M(k)$ is satisfied. Similar remark holds for the min value as well.

*Lemma 4.3:* Consider the update rule (1) and (2) with the assumptions of Theorem 3.1 with initial time $k$ such that $m(k) = \min_{j \in V} \frac{x_j(k)}{y_j(k)} < \max_{j \in V} \frac{x_j(k)}{y_j(k)} = M(k)$. Then for all $k' \geq k + D(1+\bar{\tau})$ and for all $i \in V$,

$$\frac{x_i(k')}{y_i(k')} < M(k), \text{ and,} \quad (8)$$

$$\frac{x_i(k')}{y_i(k')} > m(k). \quad (9)$$

*Proof:* Consider a node $j \in V$. There exists a node $i \in V$ such that $\frac{x_i(k)}{y_i(k)} < M(k)$. Since the graph is strongly connected, there exists a directed path from node $i$ to node $j$ in $G$. Let this path be denoted as $(m_1, i), (m_2, m_1), ..., (j, m_{d_j-1})$. Then,

$$\frac{x_{m_1}(k+\bar{\tau}+1)}{y_{m_1}(k+\bar{\tau}+1)} = \frac{p_{m_1 i} \sum_{r=0}^{\bar{\tau}} I_{k+\bar{\tau}-r, m_1 i}(r) x_i(k+\bar{\tau}-r) + \sum_{l \in N_{m_1}} p_{m_1 l} \sum_{r=0}^{\bar{\tau}} I_{k+\bar{\tau}-r, m_1 l}(r) x_l(k+\bar{\tau}-r)}{p_{m_1 i} \sum_{r=0}^{\bar{\tau}} I_{k+\bar{\tau}-r, m_1 i}(r) y_i(k+\bar{\tau}-r) + \sum_{l \in N_{m_1}} p_{m_1 l} \sum_{r=0}^{\bar{\tau}} I_{k+\bar{\tau}-r, m_1 l}(r) y_l(k+\bar{\tau}-r)} \quad (10)$$

Since $k \leq k + \bar{\tau} - r \leq k + \bar{\tau}$, using *Lemma 4.2*, it follows that $\frac{x_i(k+\bar{\tau}-r)}{y_i(k+\bar{\tau}-r)} < M(k)$. Then it follows that

$$x_i(k + \bar{\tau} - r) < M(k) y_i(k + \bar{\tau} - r). \tag{11}$$

Also, from *Lemma 4.1*, it follows that $\frac{x_l(k+\bar{\tau}-r)}{y_l(k+\bar{\tau}-r)} \leq M(k)$. Then it follows that

$$x_l(k + \bar{\tau} - r) < M(k) y_l(k + \bar{\tau} - r). \tag{12}$$

Using (11) and (12) in (10), it follows that

$$\frac{x_{m_1}(k+\bar{\tau}+1)}{y_{m_1}(k+\bar{\tau}+1)} < \frac{p_{m_1 i} \sum_{r=0}^{\bar{\tau}} I_{k+\bar{\tau}-r, m_1 i}(r) M(k) y_i(k+\bar{\tau}-r) + \sum_{l \in N_{m_1}} p_{m_1 l} \sum_{r=0}^{\bar{\tau}} I_{k+\bar{\tau}-r, m_1 l}(r) M(k) y_l(k+\bar{\tau}-r)}{p_{m_1 i} \sum_{r=0}^{\bar{\tau}} I_{k+\bar{\tau}-r, m_1 i}(r) y_i(k+\bar{\tau}-r) + \sum_{l \in N_{m_1}} p_{m_1 l} \sum_{r=0}^{\bar{\tau}} I_{k+\bar{\tau}-r, m_1 l}(r) y_l(k+\bar{\tau}-r)}$$

and hence, $\frac{x_{m_1}(k+\bar{\tau}+1)}{y_{m_1}(k+\bar{\tau}+1)} < M(k)$. Using *Lemma 4.2*, it further follows that for all $k' \geq k + \bar{\tau} + 1$, $\frac{x_{m_1}(k')}{y_{m_1}(k')} < M(k)$.

Similarly, it follows that for all $k' \geq k + 2(\bar{\tau} + 1)$, $\frac{x_{m_2}(k')}{y_{m_2}(k')} < M(k)$ and that for all $k' \geq k + d_j(\bar{\tau} + 1)$, $\frac{x_j(k')}{y_j(k')} < M(k)$. Note that since $D \geq d_j$, it follows that $D(\bar{\tau} + 1) \geq d_j(\bar{\tau} + 1)$ and hence, $k + D(\bar{\tau} + 1) \geq k + d_j(\bar{\tau} + 1)$. The condition $k' \geq k + D(\bar{\tau} + 1)$ is independent of the index $j$ and where the node $j$ was chosen arbitrarily. Thus it can be concluded that $\frac{x_i(k')}{y_i(k')} < M(k)$ for all $k' > k + D(\bar{\tau} + 1)$ and for all $i \in V$. This completes the proof. The other inequality can be proven similarly and is left to the reader. ∎

The following theorem shows that $M((l + 1)D(1 + \bar{\tau}) + (l + 1)\bar{\tau})$ (1) is a strictly decreasing sequence as a function of the index $l$. Similarly, $m((l + 1)D(1 + \bar{\tau}) + (l + 1)\bar{\tau})$ under (2) is a strictly increasing sequence as a function of the index $l$.

*Theorem 4.1:* Consider the update rule (1) and (2) with the assumptions of Theorem 3.1 and the initial ratio vector being $\frac{x(lD(1+\bar{\tau})+l\bar{\tau})}{y(lD(1+\bar{\tau})+l\bar{\tau})} := \{\frac{x_1(lD(1+\bar{\tau})+l\bar{\tau})}{y_1(lD(1+\bar{\tau})+l\bar{\tau})}, ...., \frac{x_N(lD(1+\bar{\tau})+l\bar{\tau})}{y_N(lD(1+\bar{\tau})+l\bar{\tau})}\}$ such that $\min \frac{x(lD(1+\bar{\tau})+l\bar{\tau})}{y(lD(1+\bar{\tau})+l\bar{\tau})} < \max \frac{x(lD(1+\bar{\tau})+l\bar{\tau})}{y(lD(1+\bar{\tau})+l\bar{\tau})}$, where, $l = 0, 1, 2, ....$ Then, $M((l+1)D(1+\bar{\tau})+(l+1)\bar{\tau}) < M(lD(1+\bar{\tau})+l\bar{\tau})$ and $m((l+1)D(1+\bar{\tau})+(l+1)\bar{\tau}) > m(lD(1+\bar{\tau})+l\bar{\tau})$.

*Proof:* Using $k = lD(1+\bar{\tau}) + l\bar{\tau}$ in *Lemma 4.3*, it follows that for for all $j \in V$ $\frac{x_j(lD(1+\bar{\tau})+l\bar{\tau}+D(1+\bar{\tau}))}{y_j(lD(1+\bar{\tau})+l\bar{\tau}+D(1+\bar{\tau}))} < M(lD(1+\bar{\tau})+l\bar{\tau})$. Using *Lemma 4.3*, it also follows that for for all $j \in V$ and for $p = 0, 1, 2, ..., \bar{\tau}$, $\frac{x_j((l+1)D(1+\bar{\tau})+l\bar{\tau}+p)}{y_j((l+1)D(1+\bar{\tau})+l\bar{\tau}+p)} < M(lD(1+\bar{\tau})+l\bar{\tau})$. Thus, it follows that $\max_{j \in V} \frac{x_j((l+1)D(1+\bar{\tau})+(l+1)\bar{\tau}-r)}{y_j((l+1)D(1+\bar{\tau})+(l+1)\bar{\tau}-r)} < M(lD(1+\bar{\tau})+l\bar{\tau})$. Thus, $M((l+1)D(1+\bar{\tau})+(l+1)\bar{\tau}) < M(lD(1+\bar{\tau})+l\bar{\tau})$. This completes the proof of the first inequality. The other inequality can be proved similarly and is left to the reader. ∎

*Theorem 4.2:* Consider the update rule (1) and (2) with the assumptions of Theorem 3.1. Then, $\lim_{l \to \infty} M(lD(1+\bar{\tau}) + l\bar{\tau}) = \frac{\sum_{j=1}^N x_j(0)}{N}$ and $\lim_{l \to \infty} m(lD(1+\bar{\tau}) + l\bar{\tau}) = \frac{\sum_{j=1}^N x_j(0)}{N}$.

*Proof:* It follows from Theorem 3.1 that, $\lim_{k \to \infty} \frac{x_j(k)}{y_j(k)} = \frac{\sum_{j=1}^N x_j(0)}{N}$ for all $j \in V$. This implies that for any $j \in V$ and for any given $\epsilon > 0$, there exists $K$ such that for all $k \geq K$, $|\frac{x_j(k)}{y_j(k)} - \frac{\sum_{j=1}^N x_j(0)}{N}| < \epsilon$. This implies that for all $j \in V$, $-\epsilon < \frac{x_j(k)}{y_j(k)} - \frac{\sum_{j=1}^N x_j(0)}{N} < \epsilon$. It follows that for $r \in \{0, 1, 2, ..., \bar{\tau}\}$ and for any $j \in V$, $-\epsilon < \frac{x_j(k-r)}{y_j(k-r)} - \frac{\sum_{j=1}^N x_j(0)}{N} < \epsilon$. This implies that $-\epsilon < \max_{\substack{j \in V \\ r = \{0,1,2,...,\bar{\tau}\}}} \frac{x_j(k-r)}{y_j(k-r)} - \frac{\sum_{j=1}^N x_j(0)}{N} < \epsilon$. Thus, there exists $K$ such that $k \geq K$ implies that $|M(k) - \frac{\sum_{j=1}^N x_j(0)}{N}| < \epsilon$. Similarly, $|m(k) - \frac{\sum_{j=1}^N x_j(0)}{N}| < \epsilon$. This implies that $\lim_{k \to \infty} M(k) = \frac{\sum_{j=1}^N x_j(0)}{N}$ and $\lim_{k \to \infty} m(k) = \frac{\sum_{j=1}^N x_j(0)}{N}$. Since, $\{M(lD(1+\bar{\tau})+l\bar{\tau})\}_{l=0,1,2,..}$ and $\{m(lD(1+\bar{\tau})+l\bar{\tau})\}_{l=0,1,2,..}$ are sub-sequences of convergent sequences, they converge to the same limit. Hence, it follows that $\lim_{l \to \infty} M(lD(1+\bar{\tau})+l\bar{\tau}) = \frac{\sum_{j=1}^N x_j(0)}{N}$ and $\lim_{l \to \infty} m(lD(1+\bar{\tau})+l\bar{\tau}) = \frac{\sum_{j=1}^N x_j(0)}{N}$.

*Corollary 4.1:* Consider the update rule (1) and (2) with the assumptions of Theorem 3.1. Then, $\lim_{l \to \infty} M(lD(1+\bar{\tau}) + l\bar{\tau}) - \lim_{l \to \infty} m(lD(1+\bar{\tau}) + l\bar{\tau}) = 0$.

*Proof:* The proof directly follows from *Theorem 4.2*. ∎

The following subsection introduces Maximum/ Minimum consensus protocols and the preceding results are used to design a finite time stopping criterion for ratio consensus given by (1) and (2).

## A. Maximum Consensus Protocol

The Maximum Consensus Protocol denoted by MXP computes the maximum of the given initial node conditions $z(0) = [z_1(0)\ z_2(0)\ ...\ z_N(0)]^T$ in a distributed manner. It takes $z(0)$ as an input and generates a sequence of node values based on the following update rule for node $i$,

$$z_i(k\bar{\tau} + l) = z_i(k\bar{\tau} + l - 1), l \in \{k+1, \cdots, k+\bar{\tau}\}, \tag{13}$$

$$z_i((k+1)(\bar{\tau}+1)) =$$
$$\max_{j \in N_i \cup \{i\}} \{z_j((k+1)(\bar{\tau}+1) - (r+1))I_{k-r,ij}(r)\}_{r=0,1,...,\bar{\tau}},\ k \geq 0. \tag{14}$$

Note that (13) maintains value of $z_m$ at $z_m((k-1)\bar{\tau})$ till the $k^{th}$ epoch $k\bar{\tau} + l, l \in \{1,2,...,\bar{\tau}\}$ ends. On the other hand (14) updates $z_m$ at time instances which are multiples of $\bar{\tau} + 1$ based on recent information from the neighbors and itself. Effectively every $z_m$ update takes place once after every $\bar{\tau}$ iterations.

*Lemma 4.4:* Consider a graph $G = \{V, E\}$ with time-varying delays with uniform bound $\bar{\tau}$ and an update rule for the Maximum Consensus Protocol *(MXP)* given by *(13)* and *(14)*. Let $\tilde{z} := \max_{i \in V} z_i(0)$.

1) Then, for all $k' > 0, \max_{i \in V} z_i(k') = \tilde{z}$.
2) Let for some $k$, $z_j(k) = \tilde{z}$, that is, node $j$ has the maximum value in the network at the $k^{th}$ time instant. Then, for all instants $k' > k, z_j(k') = \tilde{z}$, that is node $j$ continues to be the maximum for $k' > k$.

*Proof:* (a) The proof is left to the reader. (b) The proof follows from the fact that if the maximum value at the current iteration is held at node $j$, then node $j$ continues to hold the maximum value in future iterations as well, as the update step of the MXP (13) and (14) includes the past value of node $j$. ∎

*Lemma 4.5:* Consider an graph $G = \{V, E\}$ with possibly time-varying delays with uniform bound $\bar{\tau}$ and an update rule for the Maximum Consensus Protocol *(MXP)* given by *(13)* and *(14)*. Let $z_{\pi_1}(0) = \max_{j \in V} z_j(0)$. Then, for all $k \geq D(1+\bar{\tau})$, and any $m \in V$, $z_m(k) = z_{\pi_1}(0)$.

*Proof:* As $z_{\pi_1}(0) = \max_{j \in V} z_j(0)$ it follows from *Lemma 4.4* that $z_{\pi_1}(k) = z_{\pi_1}(0)$ for all $k \geq 0$. Consider any node $i \in V$. Since the graph $G$ is connected, there exists a path $(\pi_1, \pi_2)(\pi_2, \pi_3)...(\pi_d, i)$ connecting $i$ and $\pi_1$. It follows from the update rule (13) and (14) that within $\bar{\tau}$ iterations, $\pi_2$ will have received the value $z_{\pi_1}(0)$ and thus, $z_{\pi_2}(\bar{\tau}+1) = z_{\pi_1}(0)$; and for any $k \geq \bar{\tau}+1$, $z_{\pi_2}(k) = z_{\pi_1}(0)$. Using the above steps for $\pi_3, \pi_4, ..., \pi_d, i$, it follows that, $z_i(k) = z_{\pi_1}(0)$ for any $k \geq d(\bar{\tau}+1)$. Thus if $k \geq D(\bar{\tau}+1) \geq d(\bar{\tau}+1)$; $z_i(k) = z_{\pi_1}(0)$. Since $D$ is the diameter of the graph $G$, it follows that, for $k \geq D(\bar{\tau}+1)$, $z_m(k) = z_{\pi_1}(0) = \max_{j \in V} z_j(0)$ for all $m \in V$. This proves the theorem. ∎

## B. Minimum Consensus Protocol

The Minimum Consensus Protocol denoted by MNP computes the minimum of the given initial node conditions $w(0) = [w_1(0)\ w_2(0)\ ...\ w_N(0)]^T$ in a distributed manner. It takes $w(0)$ as an input and generates a sequence of node values based on the following update rule for node $i$,

$$w_i(k\bar{\tau}+l) = w_i(k\bar{\tau}+l-1), l \in \{k+1, \cdots, k+\bar{\tau}\}, \tag{15}$$

$$w_i((k+1)(\bar{\tau}+1)) =$$

$$\min_{j \in N_i \cup \{i\}} \{w_j((k+1)(\bar{\tau}+1) - (r+1))I_{k-r,ij}(r)\}_{r=0,1,...,\bar{\tau}},\ k \geq 0. \tag{16}$$

*Lemma 4.6:* Consider a graph $G = \{V, E\}$ with time-varying delays with uniform bound $\bar{\tau}$ and an update rule for the Minimum Consensus Protocol *(MNP)* given by *(15)* and *(16)*. Let $\tilde{w} := \min_{i \in V} w_i(0)$.

1) Then, for all $k' > 0$, $\min_{i \in V} w_i(k') = \tilde{y}$.
2) Let for some $k$, $w_j(k) = \tilde{w}$, that is, node $j$ has the minimum value in the network at the $k^{th}$ time instant. Then, for all instants $k' > k$, $w_j(k') = \tilde{w}$, that is node $j$ continues to be the minimum for $k' > k$.

*Proof:* The proof is similar to the proof of *Lemma 4.4*. ∎

*Lemma 4.7:* Consider a graph $G = \{V, E\}$ with possibly time-varying delays with uniform bound $\bar{\tau}$ and an update rule for the Minimum Consensus Protocol *(MNP)* given by *(15)* and *(16)*. Let $w_{\pi_1}(0) = \min_{j \in V} w_j(0)$. Then, for all $k \geq D(1+\bar{\tau})$, and any $m \in V$, $w_m(k) = w_{\pi_1}(0)$.

*Proof:* The proof is similar to the proof of *Lemma 4.5*. ∎

## C. Distributed Finite Time Termination Algorithm for Ratio Consensus

In this section, we propose an algorithm using the Max-Min Consensus for stopping the ratio consensus protocol in finite time based on a user specified threshold $\rho$. This algorithm first appeared in [2] as an extension of the Max-Min consensus based finite time termination of averaging consensus [3]; however, analytical guarantees were not provided. Moreover, the algorithm presented in [2] and [3] did not consider delays present on the communication channels. In [13], the authors have proposed a finite-time algorithm for distributed averaging using Max-Min Consensus in presence of fixed delays on the communication channels. Here, we deal with the general case when possibly time varying delays can be present on the communication channels. The initial conditions for the MXP and MNP protocols are set as the initial ratios held by the nodes. The MXP and MNP protocols are re-initialized at $k = nD(1+\bar{\tau})+n\bar{\tau}$, where $n = 1, 2, ...$, with $z(nD(1+\bar{\tau})+n\bar{\tau}) = \frac{x(nD(1+\bar{\tau})+n\bar{\tau})}{y(nD(1+\bar{\tau})+n\bar{\tau})}$ and $w(nD(1+\bar{\tau})+n\bar{\tau}) = \frac{x(nD(1+\bar{\tau})+n\bar{\tau})}{y(nD(1+\bar{\tau})+n\bar{\tau})}$ respectively. Define $\bar{\alpha}_i(nD(1+\bar{\tau})+n\bar{\tau}) := \max z(nD(1+\bar{\tau})+n\bar{\tau})$, $\underline{\alpha}_i(nD(1+\bar{\tau})+n\bar{\tau}) = \min w(nD(1+\bar{\tau})+n\bar{\tau})$ and $\beta_i(nD(1+\bar{\tau})+n\bar{\tau}) := \bar{\alpha}_i - \underline{\alpha}_i$.

*Theorem 4.3: Algorithm 1* converges in finite-time.

*Proof:* It follows from *Corollary 4.1* that $M(nD(1+\bar{\tau})+n\bar{\tau}) - m(nD(1+\bar{\tau})+n\bar{\tau}) \to 0$ as $n \to \infty$. Thus, for any given $\rho > 0$, there exists an integer $n \geq t(\rho)$ such that $|M(nD(1+\bar{\tau})+n\bar{\tau}) - m(nD(1+\bar{\tau})+n\bar{\tau})| < \rho$ for all nodes in

**Algorithm 1:** Finite-time termination of consensus in presence of uniformly bounded delays

**Input:**
  $x(0), y(0) = 1$ ;                                                 // Initial condition
  $D, \bar{\tau}, \rho, P$;
**Initialize:**
  $k := 0$;
  $z_i(0) := x_i(0)/y_i(0)$;
  $w_i(0) := x_i(0)/y_i(0)$;
  $\theta := 1$;
  $l := 1$;
  $\psi = 0$;
**Repeat:**
  $x_i(k+1) := p_{ii} x_i(k) + \sum_{j \epsilon N_i} p_{ij} \sum_{r=0}^{\bar{\tau}} x_j(k-r) I_{k-r,ij}(r)$;
  $y_i(k+1) := p_{ii} y_i(k) + \sum_{j \epsilon N_i} p_{ij} \sum_{r=0}^{\bar{\tau}} y_j(k-r) I_{k-r,ij}(r)$;
  **if** $k+1 = \psi + l(1+\bar{\tau})$ **then**
    /* maximum and minimum consensus updates given by (14) and (16) for each
       node $i \in V$                                                                             */
    $z_i(k+1) := \max_{j \in N_i^- \cup \{i\}} \{z_j(k-r) I_{k-r,ij}(r)\}_{r=0,1,\ldots,\bar{\tau}}$;
    $w_i(k+1) := \min_{j \in N_i^- \cup \{i\}} \{w_j(k-r) I_{k-r,ij}(r)\}_{r=0,1,\ldots,\bar{\tau}}$;
    $l := l+1$
  **end**
  **if** $k+1 = \theta((D+1)(1+\bar{\tau}))$ **then**
    $\bar{\alpha}_i := z_i(k+1)$;
    $\underline{\alpha}_i := w_i(k+1)$;
    $\beta_i := \bar{\alpha}_i - \underline{\alpha}_i$;
    **if** $\beta_i < \rho$ **then**
      **break** ;                                                       // stop $x_i$, $y_i$ and $z_i$ updates
    **else**
      set $z_i(k+1) := x_i(k+1)/y_i(k+1)$;
      $w_i(k+1) := x_i(k+1)/y_i(k+1)$;
      $\theta = \theta + 1$;
      $\psi = \psi + k$;
      $l := 1$ ;                                                        // Reset
    **end**
  **end**
  **broadcast:** $x_i(k+1), y_i(k+1), z_i(k+1), w_i(k+1)$
  $k := k+1$

the graph $G = \{V, E\}$.

∎

*Remark 3:* It has been shown in [16] that no finite time distributed averaging protocol is possible without the knowledge of some global network parameter if each node only knows its history of estimates. The only global parameters needed for this algorithm is the knowledge of graph diameter and an upper bound on the maximum delay in the network. However, it should be noted that each node does not need to know the actual diameter but some upper bound of the diameter. One of the many algorithms in literature which can accomplish the task of finding diameters is shown in [17]. Each node can first run this protocol to determine the diameter in a distributed manner.

## V. SIMULATION RESULTS

In this section we present simulation and experimental results validating the algorithm in the previous section.

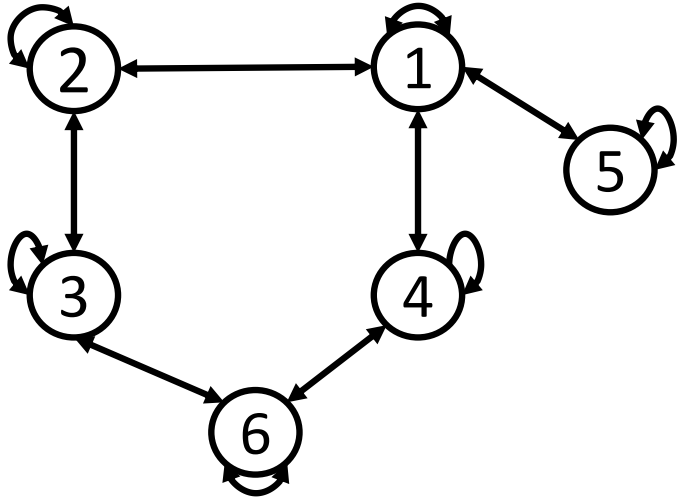

Fig. 1: A 6 node network.

*A. Simulation Results*

Consider the 6 node network shown in Fig. 1 with the weight matrix,

$$P = \begin{bmatrix} 1/4 & 1/3 & 0 & 1/3 & 1/2 & 0 \\ 1/4 & 1/3 & 1/3 & 0 & 0 & 0 \\ 0 & 1/3 & 1/3 & 0 & 0 & 1/3 \\ 1/4 & 0 & 0 & 1/3 & 0 & 1/3 \\ 1/4 & 0 & 0 & 0 & 1/2 & 0 \\ 0 & 0 & 1/3 & 1/3 & 0 & 1/3 \end{bmatrix},$$

and the initial condition of the numerator states, $x(0) = [0, 100, 200, 300, 400, 500]$ and that of the denominator states $y(0) = [1, 1, 1, 1, 1, 1]$. It is assumed that at each link at each time instant, the delay is an integer upper bounded by $\bar{\tau} = 3$. All the delays are chosen with uniform probability $1/4$ for the simulation. The application of Algorithm 1 with $\rho = 1$ and $D = 3$ results in distributed finite time termination of computations performed by the agents in 80 iterations as shown in Fig. 2. Observe that the ratio of the nodal states are close to the average of the numerator initial conditions.

## VI. CONCLUSIONS

This article presents the application of ratio consensus algorithm with max-min consensus to distributively compute average of initial conditions in finite time in the presence of time-varying delays with the possibility of multiple packets of information from same neighboring node being received at any time. The analytical results have been proven rigorously to arrive at a finite-time termination algorithm. The simulation results further establish the proposed algorithm as an effective one with considerable potential for application in practical settings such as in resource allocation in micro grid setting.

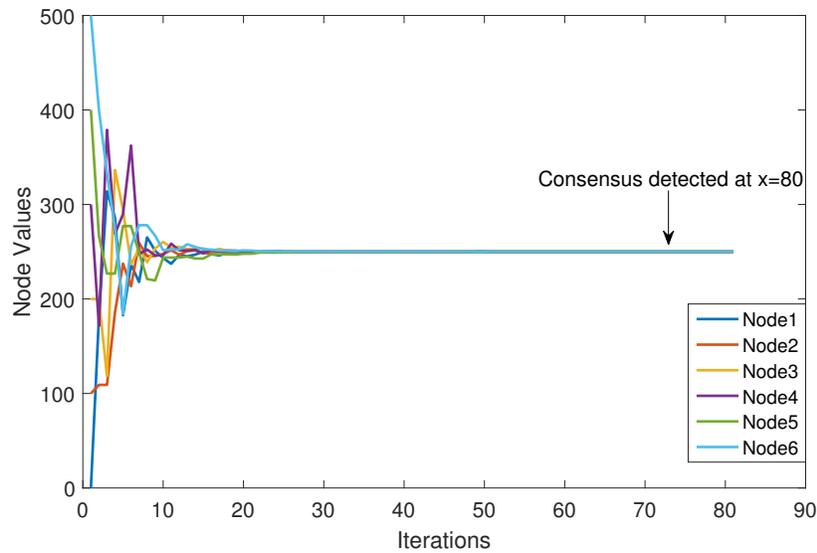

Fig. 2: Finite time termination of ratio consensus on the 6 node ring network in 80 iterations with $\rho = 1$.

## VII. Acknowledgments

The authors acknowledge the support of ARPA-E for supporting this research through the project titled 'A Robust Distributed Framework for Flexible Power Grids' via grant no. DE-AR000071 and Xcel Energy's Renewable Development Fund.